\documentclass[12pt]{article}
\usepackage[dvips]{graphicx}

\DeclareMathAlphabet{\bm}{OML}{cmm}{b}{it}
\newcommand{\startproof}{\noindent\textit{Proof.} }
\newcommand{\qed}{\ifmmode \eqno\rule{1.3ex}{1.3ex}
                  \else\hspace*{\fill}\rule{1.3ex}{1.3ex}
                  \fi}
\newtheorem{proposition}{Proposition}
\newtheorem{theorem}[proposition]{Theorem}

\newtheorem{remark}[proposition]{Remark}

\begin{document}
\title{Decoding method for generalized algebraic
geometry codes\thanks{%
This paper is registered at the LANL eprint server
\texttt{http://\allowbreak arXiv.\allowbreak org/}.}}
\author{%
\begin{tabular}{c@{\hspace*{1cm}}c}Ryutaroh Matsumoto&
Masakuni Oishi\thanks{All correspondence should be sent to the first author.}\\
\texttt{ryutaroh@rmatsumoto.org}&\end{tabular}\\
Department of Communications and Integrated Systems\\
Tokyo Institute of Technology, 152-8552 Japan}
\date{April 24, 2001}
\maketitle
\begin{abstract}
We propose a decoding method for the generalized algebraic geometry
codes proposed by Xing et~al.
To show its practical usefulness,
we give an example of generalized algebraic geometry codes
of length $567$ over $\mathbf{F}_8$ whose numbers of correctable
errors by the proposed method
are larger than the shortened codes of the
primitive BCH codes of length $4095$ in the most range of dimension.
\end{abstract}
\section{Introduction}
Recently Xing, Niederreiter and Lam \cite{xing99}
introduced
a construction of linear codes from
algebraic function fields and places of degrees larger than one.
They call their new construction \emph{generalized algebraic geometry
codes}, which is a generalization of
functional algebraic geometry codes
first explicitly defined in \cite{vanlint87,stichtenoth88}.
The lower bound for the minimum distance by Xing et~al.\ was
improved by \"Ozbudak and Stichtenoth \cite{ozbudak99}.

Nobody has proposed a decoding method that corrects errors
up to half the designed minimum distance by \"Ozbudak and
Stichtenoth. We give a decoding method
that corrects errors almost half the designed minimum
distance.

Although several examples were given \cite{xing99},
the usefulness of generalized algebraic geometry codes
has not been clearly demonstrated.
We give an example of
generalized algebraic geometry codes whose numbers of correctable
errors by the proposed method are larger than the shortened codes of BCH codes
of the same dimension (see Fig.\ \ref{graph}).

It should be noted that Kaida et~al.\ \cite{kaida95}
also considered a code construction using places of
degree larger than one,
which is a special case of generalized algebraic geometry codes.
However, they gave neither lower bound for the minimum distance
nor interesting example.

\section{Notations}
Let us fix notations used in this paper.
We shall always consider linear codes over the finite field
$\mathbf{F}_q$ with $q$ elements.
Notations follow those in Stichtenoth's textbook \cite{bn:stichtenoth}.
Let $F/\mathbf{F}_q$ be an algebraic function field of one variable
over the full constant field $\mathbf{F}_q$.
Let $P_1$, \ldots, $P_s$ be pairwise distinct places of $F/\mathbf{F}_q$,
not necessarily of degree one.
Let $G$ be a divisor of $F/\mathbf{F}_q$ whose support contains
none of $P_1$, \ldots, $P_s$.

For ease of presentation
we shall define a subclass of generalized algebraic geometry codes
in a slightly different manner from the original \cite{xing99}.
Let $m$ be the least common multiple of $\deg P_1$, \ldots, $\deg P_s$.
Let $\tilde{F} = F\mathbf{F}_{q^m}$, that is, the constant field
extension of $F$.
Let $\tilde{P}_{i,1}$, \ldots, $\tilde{P}_{i,\deg P_i}$ be the extensions
of $P_i$ in $\tilde{F}/F$.
The residue class field $O_{P_i}/P_i$ can be regarded as a subfield
of $O_{\tilde{P}_i,1}/ \tilde{P}_{i,1}$ via the $\mathbf{F}_q$-embedding
$\iota_i$
sending $(x \bmod P_i)\in O_{P_i}/P_i$ to $(x \bmod \tilde{P}_{i,1})
\in O_{\tilde{P}_i,1}/ \tilde{P}_{i,1}$.
Let $\iota_i(F_{P_i})$ be the image of $F_{P_i} = O_{P_i}/P_i$
under $\iota_i$. We have the field isomorphism
\begin{equation}
\iota_i(F_{P_i}) \simeq \mathbf{F}_{q^{\deg P_i}}. \label{iota}
\end{equation}

Let $\pi_i$ be an $\mathbf{F}_q$-linear isomorphism from
$\iota_i(F_{P_i})$ to $\mathbf{F}_q^{\deg P_i}$.
We shall consider the generalized algebraic geometry code defined by
\[
C = \{ (\pi_1(f(\tilde{P}_{1,1})), \pi_2(f(\tilde{P}_{2,1})),
\ldots, \pi_s(f(\tilde{P}_{s,1}))) \mid
f \in \mathcal{L}(G)\},
\]
where $\mathcal{L}(G)$ is $\{ f \in F \mid$ the principal divisor of
$f \geq -G\}$.
The length of $C$ is given by
$n = \deg P_1 + \cdots + \deg P_s$.
The dimension of $C$ is given by
$\dim G - \dim (G-(P_1 + \cdots + P_s))$.

In the original definition of generalized algebraic geometry codes
\cite{xing99},
$\pi_i$ is defined as an injective $\mathbf{F}_q$-linear map from
$F_{P_i}$ to $\mathbf{F}_q^{\gamma_i}$, where $\gamma_i$ is an integer
$\geq \deg P_i$. Our definition is equivalent to
the case of $\gamma_i = \deg P_i$.

\section{Decoding method}
In this section we propose a decoding method for $C$
using a decoding method for a code over $\mathbf{F}_{q^m}$.

The field $\iota_i(F_{P_i})$ is isomorphic to
$\mathbf{F}_q^{\deg P_i}$ as an $\mathbf{F}_q$-space.
We shall consider the $\mathbf{F}_q$-space
\[
\tilde{C} = \{ (f(\tilde{P}_{1,1}),
f(\tilde{P}_{2,1}), \ldots, f(\tilde{P}_{s,1})) \mid
f \in \mathcal{L}(G) \}
\]
that can be regarded as an $\mathbf{F}_q$-subspace of
\[
\Lambda = \mathbf{F}_{q^{\deg P_1}} \times \mathbf{F}_{q^{\deg P_2}}
\times \cdots \times \mathbf{F}_{q^{\deg P_s}}
\]
by isomorphisms (\ref{iota}).

For an array $\bm{x} = (x_1$, \ldots, $x_s)$, where $x_i$ is an element
in some field,
we define the Hamming weight of $\bm{x}$ to be the number of nonzero
component in $\bm{x}$.

Let $\sigma$ be the Frobenius automorphism of $\mathbf{F}_{q^m}$
sending $\alpha$ to $\alpha^q$. $\sigma$ can be extended to the automorphism
of $\tilde{F}/F$ in a natural manner.
We may assume $\sigma \tilde{P}_{i,j} = \tilde{P}_{i,j+1}$
without loss of generality \cite[Theorem III.7.1]{bn:stichtenoth}.

We define the following $\mathbf{F}_q$-spaces:
\begin{eqnarray*}
\Lambda_\mathrm{ex} &=& \mathbf{F}_{q^{\deg P_1}}^{\deg P_1} \times \mathbf{F}_{q^{\deg P_2}}^{\deg P_2}
\times \cdots \times \mathbf{F}_{q^{\deg P_s}}^{\deg P_s},\\
\tilde{C}_\mathrm{ex} &=&
\{ (f(\tilde{P}_{1,1}),f(\tilde{P}_{1,2}), \ldots, f(\tilde{P}_{1,\deg P_1}), f(\tilde{P}_{2,1}),
\ldots, f(\tilde{P}_{s,\deg P_s})) \mid f \in \mathcal{L}(G) \}.
\end{eqnarray*}
By isomorphisms (\ref{iota}), $\tilde{C}_\mathrm{ex}$ can be regarded
as an $\mathbf{F}_q$-subspace of $\Lambda_\mathrm{ex}$.

For $f\in F$ with $v_{P_i}(f) \geq 0$,
we have $\sigma(f(\tilde{P}_{i,j})) = (\sigma f)(\sigma \tilde{P}_{i,j}) = f(\sigma \tilde{P}_{i,j})$
\cite[Proof of Lemma III.5.2 (c)]{bn:stichtenoth}.
So we have $(f(\tilde{P}_{1,1})$, $f(\tilde{P}_{1,2})$, \ldots, $f(\tilde{P}_{1,\deg P_1})$, 
$f(\tilde{P}_{2,1})$,
\ldots, $f(\tilde{P}_{s,\deg P_s})) = 
(f(\tilde{P}_{1,1})$, $\sigma(f(\tilde{P}_{1,1}))$, \ldots, $\sigma^{\deg P_1 - 1}(
f(\tilde{P}_{1,1}))$, 
$f(\tilde{P}_{2,1})$, \ldots, $\sigma^{\deg P_s-1}(f(\tilde{P}_{s,\deg P_s})))$.

Consider an $\mathbf{F}_q$-linear map $\varphi : \Lambda \rightarrow
\Lambda_\mathrm{ex}$ sending $(x_1$, \ldots, $x_s)$ to
$(x_1$, $\sigma x_1$, \ldots, $\sigma^{\deg P_1 - 1} x_1$, $x_2$,
$\sigma x_2$,
\ldots, $\sigma^{\deg P_s - 1}x_s)$.
Then $\varphi$ is injective and $\varphi(\tilde{C}) = \tilde{C}_\mathrm{ex}$.

Let $\mathrm{Con}_{\tilde{F}/F}(G)$ be the conorm of $G$ in $\tilde{F}/F$
\cite[Definition III.1.8]{bn:stichtenoth},
then $\mathcal{L}(G)$ is an $\mathbf{F}_q$-subspace of
$\mathcal{L}(\mathrm{Con}_{\tilde{F}/F}(G))$ \cite[Theorem III.6.3 (d)]{bn:stichtenoth}.
Thus $\tilde{C}_\mathrm{ex}$ is a subcode of the ordinary algebraic geometry
code
\[
\tilde{C}_\mathrm{ord} =
\{f(\tilde{P}_{1,1}), f(\tilde{P}_{1,2}), \ldots, f(\tilde{P}_{s,\deg P_s})) \mid
f \in \mathcal{L}(\mathrm{Con}_{\tilde{F}/F}(G)) \}.
\]

Suppose that a codeword $\bm{c} \in C$ is sent
and $\bm{r} = \bm{c} + \bm{e} \in \mathbf{F}_q^n$ is received.
Let $\pi$ be an $\mathbf{F}_q$-linear map from $\Lambda$ to
$\mathbf{F}_q^n$ defined by $\pi_1 \times \cdots \times \pi_s$.
When the number of errors is not too large,
we can find $\varphi(\pi^{-1}(\bm{c}))$ from $\varphi(\pi^{-1}(\bm{r}))$
by using a decoding algorithm for $\tilde{C}_\mathrm{ord}$.
We shall analyze the number of correctable errors
by a $t$-error correcting algorithm for $\tilde{C}_\mathrm{ord}$.

We shall relate the Hamming weight of $\bm{x}$ in $\Lambda$ and
$\varphi(\bm{x})$ in $\Lambda_\mathrm{ex}$.
We define $\nu_i = \sharp \{ j \mid \deg P_j = i \}$,
where $\sharp$ denotes the number of elements in a set,
and $\mu = \max\{ \deg P_i \mid i=1$, \ldots, $s\}$.
\begin{proposition}\label{lower}
Suppose that
\begin{equation}
\sum_{i=a}^\mu i\nu_i \leq w(\varphi(\bm{x})) < \sum_{i=a+1}^\mu i\nu_i,
\label{adefinition}
\end{equation}
where $\bm{x} \in \Lambda$, and $w(\cdot)$ denotes the Hamming
weight of $\cdot$.
Then we have
\begin{eqnarray}
w(\bm{x}) &\geq & \min\{ \sharp S \mid 
S \subseteq \{1,\ldots,s\},\;
\sum_{i\in S} \deg P_i \geq w(\varphi(\bm{x}))\} \label{weightrelation}\\
&=& \left\lceil \frac{w(\varphi(\bm{x})) - \sum_{i=a+1}^\mu (i-a)\nu_i}{a}\right\rceil.\nonumber
\end{eqnarray}
\end{proposition}
\startproof
Let $\bm{x} = (x_1$, \ldots, $x_s)$.
If $x_i \neq 0$ then
$\sigma^j x_i \neq 0$ for any $j$.
It follows that
\[
w(\varphi(\bm{x})) = \sum_{x_i \neq 0} \deg P_i.
\]
Thus we have
\begin{eqnarray*}
w(\bm{x}) &\geq& \min\{ \sharp S \mid S \subseteq \{1,\ldots,s\},\;
\sum_{i \in S} \deg P_i = w(\varphi(\bm{x}))\} \\
&\geq& \min\{ \sharp S \mid S \subseteq \{1,\ldots,s\},\;
\sum_{i \in S} \deg P_i \geq w(\varphi(\bm{x}))\}.
\end{eqnarray*}

We may assume $\deg P_1 \geq \cdots \geq \deg P_s$
without loss of generality.
Define $\ell$ by
\[
\sum_{i=1}^{\ell-1} \deg P_i < w(\varphi(\bm{x})) \leq
\sum_{i=1}^\ell \deg P_i.
\]
We can easily see
\[
\ell = \min\{ \sharp S \mid S \subseteq \{1,\ldots,s\},\;
\sum_{i \in S} \deg P_i \geq w(\varphi(\bm{x}))\}.
\]
By the definitions of $a$ and $\ell$,
we have $a = \deg P_\ell$, and
\begin{eqnarray*}
\ell &=& \sum_{i=a+1}^\mu \nu_i + \left\lceil
\frac{w(\varphi(\bm{x}))-\sum_{i=a+1}^\mu i \nu_i}{a}\right\rceil \\
&=& \left\lceil \frac{w(\varphi(\bm{x})) - \sum_{i=a+1}^\mu (i-a)\nu_i}{a}
\right\rceil. \qquad\mbox{\qed}
\end{eqnarray*}

\begin{theorem}
When we have a $t$-error correcting algorithm for
$\tilde{C}_\mathrm{ord}$,
we can correct up to
\[
\left\lceil \frac{t+1 - \sum_{i=a+1}^\mu (i-a)\nu_i}{a}\right\rceil-1
\]
errors of $C$,
where $a$ is defined by
\[
\sum_{i=a}^\mu i\nu_i \leq t+1 < \sum_{i=a+1}^\mu i\nu_i.
\]
\end{theorem}
\startproof
Suppose that we sent $\bm{c} \in C$ and received $\bm{r} = \bm{c}
+\bm{e} \in \mathbf{F}_q^n$. By using the $t$-error correcting
algorithm for $\tilde{C}_\mathrm{ord}$,
we try to find $\varphi(\pi^{-1}(\bm{c}))$ from $\varphi(\pi^{-1}(\bm{r}))$,
and compute $\bm{c}$ from $\varphi(\pi^{-1}(\bm{c}))$.
In order for this method to work,
it is sufficient that $w(\varphi(\pi^{-1}(\bm{e}))) \leq t$.
We shall show that if $w(\bm{e}) \leq
\lceil (t+1 - \sum_{i=a+1}^\mu (i-a)\nu_i)/a\rceil -1$ then
$w(\varphi(\pi^{-1}(\bm{e}))) \leq t$.

We have $w(\bm{e}) \geq w(\pi^{-1}(\bm{e}))$.
Suppose that $w(\varphi(\pi^{-1}(\bm{e}))) = t+1$.
Then by Proposition \ref{lower}
\begin{equation}
w(\bm{e}) \geq \left\lceil \frac{t+1 - \sum_{i=a+1}^\mu (i-a)\nu_i}{a}\right\rceil\label{minweight},
\end{equation}
which is a contradiction.
Suppose that $w(\varphi(\pi^{-1}(\bm{e}))) > t+1$.
Then the right hand side of Eq.\ (\ref{minweight}) increases and
the same contradiction is deduced. \qed

\begin{remark}
Decoding algorithms for the functional algebraic geometry code
$\tilde{C}_\mathrm{ord}$ are proposed in
\textup{\cite{duursma93,farran00,guruswami99,ldecodepaper}.}
In \textup{\cite{ldecodepaper}} the Feng-Rao decoding algorithm
\textup{\cite{fengrao93}}
is modified for functional algebraic geometry codes.
Note that we have to represent $\tilde{C}_\mathrm{ord}$
as a residue algebraic geometry code before applying the algorithms
\textup{\cite{duursma93,farran00},} and the algorithms
\textup{\cite{guruswami99,ldecodepaper}}
are applicable only for one-point codes.
\end{remark}

\begin{remark}\label{closedformula}
\"Ozbudak and Stichtenoth \textup{\cite{ozbudak99}} showed that
the minimum distance of $\tilde{C}$ is not less than
\begin{equation}
\min\{ \sharp S \mid S \subseteq \{1,\ldots,s\},\;
\sum_{i\in S} \deg P_i \geq n-\deg G\}. \label{ozbudak}
\end{equation}
By an argument similar to the proof of Proposition \textup{\ref{lower},}
one can prove that the lower bound \textup{(\ref{ozbudak})} is
equal to
\[
\left\lceil \frac{n-\deg G - \sum_{i=a+1}^\mu (i-a)\nu_i}{a}\right\rceil,
\]
where $a$ is defined by
\[
\sum_{i=a}^\mu i\nu_i \leq n-\deg G < \sum_{i=a+1}^\mu i\nu_i.
\]
\end{remark}

\begin{remark}\label{smallestdegrees}
In order to make the number of correctable errors and
the lower bound \textup{(\ref{ozbudak})} larger,
we have to make the right hand side (RHS) of Eq.\ \textup{(\ref{weightrelation})}
larger.
The RHS of Eq.\ \textup{(\ref{weightrelation})} takes the maximum value
for a fixed code length $n$ and $w(\varphi(\bm{x}))$ when
the places $P_1$, \ldots, $P_s$ are of the smallest degrees.
\end{remark}

\section{Examples}
\subsection{Good example}
\begin{figure}[bth!]
\includegraphics[width=\linewidth]{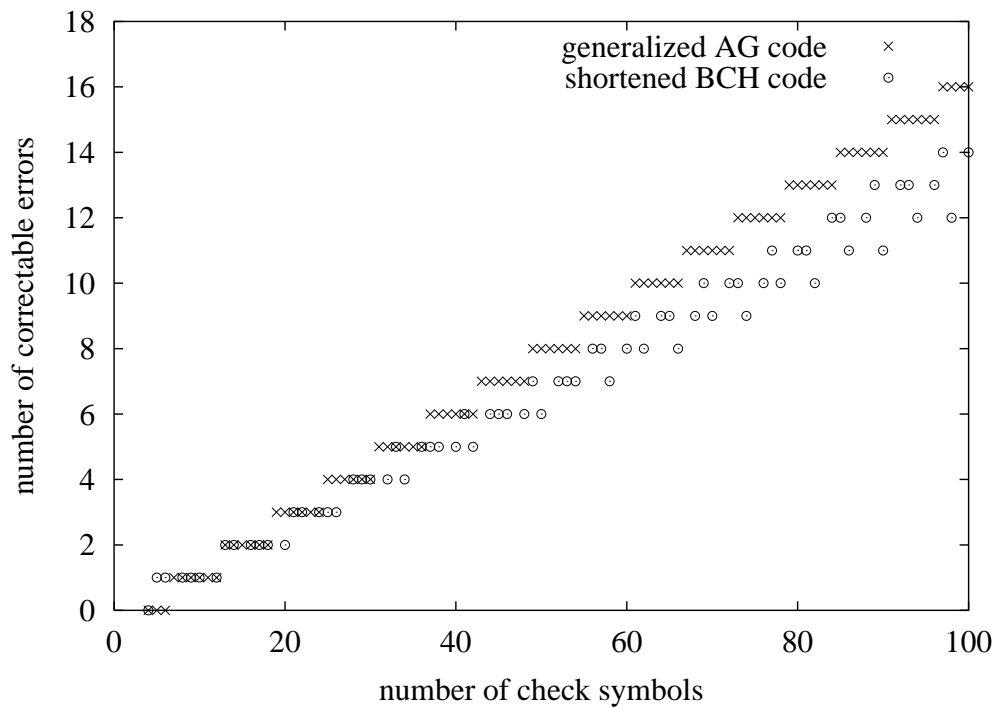}
\caption{Comparison of generalized AG codes and BCH codes}\label{graph}
\end{figure}

In this subsection we compare generalized algebraic geometry codes
and BCH codes of length $567$ over $\mathbf{F}_8$.
We construct codes from the rational function field
$\mathbf{F}_8(x)/\mathbf{F}_8$.
We take $D$ as the sum of 
$7$ places of degree $1$,
$28$ places of degree $2$,
and $168$ places of degree $3$.
We compare these generalized algebraic geometry codes
with shortened codes of primitive BCH codes of length $4095$.
For each number of check symbols,
we take a BCH code that has the largest BCH bound.
The number of check symbols and the number of correctable
errors of these codes are plotted in Fig.\ \ref{graph}.

\subsection{Bad example}
In this subsection
we give an example of generalized algebraic geometry codes
with which we cannot correct errors up to half the
designed minimum distance.
Consider the rational function field
$\mathbf{F}_{17}(x)/\mathbf{F}_{17}$.
Take $D$ as the sum of $17$ places of degree $1$
and a place of degree $4$, and
$G$ as a divisor of degree $13$.
Then the designed minimum distance is
$5$, while the number of correctable errors is $0$.

\section*{Acknowledgment}
We would like to thank Prof.\ Tomohiko Uyematsu for helpful comments.


\end{document}